\documentclass[journal,twoside,web]{ieeecolor}

\usepackage{generic}
\usepackage{cite}
\usepackage{amsmath,amssymb,amsfonts}
\usepackage{algorithmic}
\usepackage{graphicx}
\usepackage[dvipsnames]{xcolor}
\usepackage{algorithm,algorithmic}
\usepackage{hyperref}
\usepackage{textcomp}
\def\BibTeX{{\rm B\kern-.05em{\sc i\kern-.025em b}\kern-.08em
    T\kern-.1667em\lower.7ex\hbox{E}\kern-.125emX}}
\usepackage{mathtools,gensymb,mathrsfs,bm}
\usepackage{diagbox} 
\usepackage{ booktabs}
\usepackage{tikz}
\usetikzlibrary{calc}

\newcommand{\tikzmark}[1]{\tikz[overlay,remember picture] \node (#1) {};}
\newcommand{\DrawBox}[3][]{%
    \tikz[overlay,remember picture]{
    \draw[black,#1]
      ($(#2)+(-0.5em,2.0ex)$) rectangle
      ($(#3)+(0.75em,-0.75ex)$);}
}
\usepackage{etoolbox}
\makeatletter
\patchcmd{\@makecaption}
  {\scshape}
  {}
  {}
  {}
\makeatletter
\patchcmd{\@makecaption}
  {\\}
  {.\ }
  {}
  {}
\makeatother

\AtBeginEnvironment{pmatrix}{\setlength{\arraycolsep}{1.5pt}}

\newtheorem{thm}{Theorem}

\newtheorem{Lem}{Lemma}

\newcommand{\ie}{{i.e. }}

\newcommand{\e}{\epsilon}

\begin{document}
\title{\color{black}Exponential Stability and Design of Sensor Feedback Amplifiers for Fast Stabilization of Magnetizable Piezoelectric Beam Equations}
\author{Ahmet \"Ozkan \"Ozer$^{1}$, \IEEEmembership{Member, IEEE},  Ahmet Kaan Ayd\i n$^{2}$ \IEEEmembership{Member, IEEE}, and Rafi Emran$^{3}$
	\thanks{$^{1}$Department of Mathematics, Western Kentucky University,
		Bowling Green, KY 42101, USA.
		{\tt\small ozkan.ozer@wku.edu } }
	\thanks{$^{2}$Department of Mathematics, University of Maryland, Baltimore County, Baltimore, MD 21250, USA.
		{\tt\footnotesize aaydin1@umbc.edu} }
	\thanks{$^{3}$Department of Mathematical Sciences, New Jersey Institute of Technology, Newark, NJ 07102, USA.
		{\tt\footnotesize ms3532@njit.edu} }
	 }

\maketitle

\begin{abstract}
\color{black}
The dynamic partial differential equation (PDE) model governing longitudinal oscillations in magnetizable piezoelectric beams exhibits exponentially stable solutions when subjected to two boundary state feedback controllers. An analytically established exponential decay rate by the Lyapunov approach ensures stabilization of the system to equilibrium, though the actual decay rate could potentially be improved. The decay rate of the closed-loop system is highly sensitive to the choice of material parameters and the design of the state feedback amplifiers. This paper focuses on investigating the design of state feedback amplifiers to achieve a maximal exponential decay rate, which is essential for effectively suppressing oscillations in these beams. Through this design process, we explicitly determine the safe intervals of feedback amplifiers that ensure the theoretically found maximal decay rate, with the potential for even better rates. Our numerical results reaffirm the robustness of the decay rate within the chosen range of feedback amplifiers, while deviations from this range significantly impact the decay rate. To underscore the validity of our results, we present various numerical experiments.
\end{abstract}


\section{Introduction}
\label{sec:introduction}
{ \IEEEPARstart{P}iezoelectric ceramics, including lead-free materials like barium titanate, sodium potassium niobate, and sodium bismuth titanate, are renowned multifunctional smart materials that generate electric displacement in response to mechanical stress~\cite{Smith}. Their small size and high power density make them ideal for various industrial applications, including implantable biomedical devices~\cite{Dag1, Shi}, wearable interfaces with PVDF sensors~\cite{Dong}, biocompatible sensors~\cite{Biomim}, and ultrasound imagers and cleaners~\cite{Sci-D}. Their fast response, large mechanical force, and fine resolution contribute to their effectiveness~\cite{Gu}.

Consider a piezoelectric beam, clamped on one side and free on the other, with sensors for tip velocity and current. he beam is of length $L$ and thickness $h$. Assume that the transverse oscillations of the beam are negligible, so the longitudinal vibrations, in the form of expansion and compression of the center line of the beam, are the only oscillations of note together with the electromagnetic effects. Existing mathematical models often oversimplify intrinsic mechanical or electromagnetic interactions, impacting boundary feedback stabilizability~\cite{Aydin1, M-O, O-M, Sherp, Voss, Prieur}.

While electrostatic and quasi-static approaches based on Maxwell's equations are typically sufficient for non-magnetizable piezoelectric beams, assessing radiated electromagnetic power requires consideration of electromagnetic waves generated by mechanical fields~\cite{Yang}. Therefore, fully dynamic models of piezoelectric beams are essential. Denoting \( v(x,t) \) as the longitudinal oscillations and \( p(x,t) \) as the total charge accumulated at the electrodes, the equations of motion form following system of partial differential equations~\cite{M-O}
\begin{eqnarray}
\label{eq1}
\begin{array}{ll}
	\begin{array}{ll}
		\begin{bmatrix}
				\rho & 0 \\
				0 & \mu
		\end{bmatrix}
		\begin{bmatrix}
				v_{tt}  \\
				p_{tt}
		\end{bmatrix} -
		\begin{bmatrix}
				\alpha & -\gamma\beta \\
				-\gamma\beta & \beta
		\end{bmatrix}
		\begin{bmatrix}
				v_{xx}  \\
				p_{xx}
		\end{bmatrix}= \begin{bmatrix}
				0 \\
				0
			\end{bmatrix}, &\\
	\end{array}\\	
	\begin{array}{ll}
		(v,p)(0,t)=0,&  \\
		\begin{bmatrix}
			\alpha & -\gamma\beta \\
			-\gamma\beta & \beta
		\end{bmatrix}
		\begin{bmatrix}
			v_{x}  \\
			p_{x}
		\end{bmatrix} (L,t)=\begin{bmatrix}
			u_1(t) \\
			u_2(t)
		\end{bmatrix}, \quad t\in\mathbb{R}^+
	\end{array} \\
  \left[v,p,v_t,p_t\right]\left(x,0\right)=\left[v_0,p_0,v_1,p_1\right]\left(x\right), ~ x\in\left[0,L\right]
\end{array}
\end{eqnarray}
where  $\rho$, $\alpha$, $\beta$, $\gamma$, and $\mu$ are the mass density per unit volume, the elastic stiffness, the impermeability, the piezoelectric constant, and the magnetic permeability, respectively, and $u_1(t)$ and $u_2(t)$ are strain and voltage actuators, respectively.

  Define $\alpha:=\alpha_1 + \gamma^2 \beta>0$ with $\alpha_1>0$, and
  \begin{eqnarray}
	\label{constants}\begin{array}{l}
		\zeta^{\mp} = \frac{1}{\sqrt{2}}\sqrt{\frac{\alpha\mu}{\alpha_1\beta}+\frac{\rho}{\alpha_1}\mp\sqrt{\left(\frac{\alpha\mu}{\alpha_1\beta}+\frac{\rho}{\alpha_1}\right)^2-\frac{4\rho\mu}{\beta\alpha_1}}}.
	\end{array}&&
\end{eqnarray}
Note that $\sqrt{\frac{\rho}{\alpha}}$ and $\sqrt{\frac{\beta}{\mu}}$ are non-identical and represent significantly different wave propagation speeds in \eqref{eq1}.
The natural energy of the solutions is defined as
\begin{eqnarray}\label{eq4}\begin{array}{l}
	E(t)=\frac{1}{2}\int^L_0\left[\rho\left|v_t\right|^2+ \mu \left|p_t\right|^2+\alpha_1 |v_x|^2 \right.\\ \hspace{4.5cm}\left.+\beta\left|\gamma v_x-p_x\right|^2\right] dx.
\end{array}
\end{eqnarray}

The exact observability result for the model~\eqref{eq1} with $u_1(t),~u_2(t)=0$  and with only one sensor measurement is not possible \cite{M-O}. However, by employing two sensor measurements, $v_t(L,t)$ (tip velocity) and $p_t(L,t)$ (total current accumulated at electrodes), a suboptimal observation time is achieved \cite{Ramos}. This result is later refined to include the optimal observation time \cite{Wilson}.
 \begin{thm} \cite[Theorem 2.2]{Wilson} Define the Hilbert spaces $(H^1_*(0,L))=\{z\in H^1(0,L): z(0)=0\}$ and $\mathcal{H}=(H^1_*(0,L))^2\times (L^2(0,L))^2.$ For all initial data $\left[v_0,p_0,v_1,p_1\right](x)\in \mathcal{H}$, and for any $T>\frac{2L}{\max\{\zeta^-,\zeta^+\}}$,
  there exists a constant $C(T)>0$ such that the weak solutions $(v,p,v_t,p_t)\in \mathcal{H}$ of the control-free system~\eqref{eq1}, i.e., $u_1(t),~u_2(t)=0,$ satisfy
\begin{eqnarray}\begin{array}{l}
	\int^T_0\left(\rho\left|v_t(L,t)\right|^2+\mu\left|p_t(L,t)\right|^2\right)dt\geq C(T)E(0).
\end{array}
\end{eqnarray}
\end{thm}

When the observed sensor signals $v_t(L,t)$ and $p_t(L,t)$ are amplified and fed back to~\eqref{eq1}, a closed-loop system is formed. The amplification range depends on the sensor limits. While using only one sensor feedback can make the system energy dissipative, it is insufficient for exponential stabilization in $\mathcal{H}$. This is because the closed-loop system with this control design imposes strict conditions on stabilization results \cite{M-O}. Exponential stability is jeopardized for a large class of material parameters and is only achievable for a small subset \cite{MCSS}.

{
Let $\xi_1, \xi_2>0$ be the feedback amplifiers for the two sensor measurements of the closed-loop system $v_t(L,t)$ and $p_t(L,t)$, respectively. The strain and voltage actuators $u_1(t)$ and $u_2(t)$ are chosen to be proportional to the sensor feedback amplifiers,
\begin{equation}\label{c1}
    \begin{array}{ll}
		\begin{bmatrix}
			u_1(t) \\
			u_2(t)
		\end{bmatrix} =-\begin{bmatrix}
			\xi_1      &0 \\
			0       & \xi_2
		\end{bmatrix}\begin{bmatrix}
			\dot v(L,t) \\
			\dot p(L,t)
		\end{bmatrix}.
	\end{array}
\end{equation}

The exponential stability of the closed-loop system \eqref{eq1}, \eqref{c1} has been studied in~\cite{Ramos}, where the proof relies on decomposing the system into conservative and dissipative components. However, this method does not explicitly describe the decay rate or allow optimization of the feedback amplifiers.
}

{
The primary objective of this paper is to establish the existence of an exponential decay rate for the system described by \eqref{eq1} and \eqref{c1} through the meticulous construction of a Lyapunov function. Given the impracticality of traditional spectral analysis due to the strong coupling in the wave system, we adopt a multiplier approach combined with an optimization argument to define a safe range of intervals for each feedback amplifier $\xi_1$ and $\xi_2$. This ensures that the decay rate provided by the Lyapunov approach is achievable for any type of initial conditions $\mathcal{H}$.

Our methodology provides a maximal decay rate for the system, serving as an upper bound to the optimal decay rate since the actual decay rate could potentially be even better. {\color{black}Existing literature on optimal actuator designs commonly offers different approaches for infinite-dimensional systems~\cite{Ilias, Vexler2016, Mor-Dem, Privat}, and finite-dimensional systems~\cite{Munch, GeshkovskiBrunovsky}. Our approach is particularly valuable for model reductions by Finite Differences or Finite Elements~\cite{Wilson} for~\eqref{eq1}. Notably, the Lyapunov-based exponential stability proof~\cite{O-Emran2}, uniformly as the discretization parameter $h\to 0$ with the recently proposed order-reduced Finite Differences \cite{GuoWave}, closely resembles the approach employed here.}

}
\section{Exponential Stability Result}
\label{sec:ExpStab}
{For the solutions of the system~\eqref{eq1},\eqref{c1} to stabilize exponentially, the energy must be dissipative. The proof of the following dissipativity theorem is omitted.
\begin{Lem}\label{thm1}
	For all $t>0,$ the energy $E(t)$ in \eqref{eq4} is dissipative. In other words,
\[
    \frac{dE}{dt}=-\xi_1\left|\dot v(L,t)\right|^2-\xi_2 \left|\dot p(L,t)\right|^2\le 0.
\]
\end{Lem}
}
Now, let's define an energy-like functional $F(t)$ in order to establish a perturbed energy functional $E_\delta(t)$ as follows
\begin{align}\label{eq5}
		&    F(t):=\int^L_0\left(\rho v_txv_x+\mu p_txp_x\right)dx,\\
		\label{eq6}
		&       E_\delta(t):=E(t)+\delta F(t).
\end{align}
Here, $\delta>0$ will be determined as a function of sensor feedback amplifiers $\xi_1$ and $\xi_2$ later.

The following two lemmas for $F(t)$ and $E_\delta(t)$ are needed to prove our main exponential stability result. Let
 \begin{equation}\label{eta}
	\eta:=\max\left(\sqrt{\frac{\rho}{\alpha_1}}+\sqrt{\frac{\mu\gamma^2}{\alpha_1}},\sqrt{\frac{\mu}{\beta}}+\sqrt{\frac{\mu\gamma^2}{\alpha_1}}\right),
\end{equation}

\begin{Lem}\label{lem2}
	Letting $0<\delta <\frac{1}{\eta L}$, for all $t>0$, $E_\delta(t)$ in \eqref{eq6} is equivalent to $E(t)$ in \eqref{eq4}. In other words,
\begin{equation}\label{equiva}
	\left(1-\delta \eta L \right)E(t) \leq E_\delta(t) \leq\left(1+\delta \eta L \right)E(t).
\end{equation}
\end{Lem}

\begin{proof} By the H\"older's, Minkowski's, Triangle inequalities,  as well as algebraic manipulations, $F(t)$ satisfies
	\begin{equation*}
		\begin{array}{ll}
				&|F(t)|\leq L\left[\left(\int_{0}^L\rho |v_t|^2{dx}\right)^\frac{1}{2}\left(\frac{\rho}{\alpha_1}\int_{0}^L\ \alpha_1|v_x|^2{dx}\right)^\frac{1}{2}\right.\\
			& \quad \left. +\left(\int_{0}^L\mu |p_t|^2{dx}\right)^\frac{1}{2}\left\{\left(\int_0^L\mu |\gamma v_x-p_x|^2dx\right)^\frac{1}{2}\right.\right.\\
			&\qquad\left.\left.+\left(\mu\int_0^L\gamma^2|v_x|^2dx\right)^\frac{1}{2}\right\}\right]\\
			&\leq \frac{L}{2}\left[\sqrt{\frac{\rho}{\alpha_1}}\int_{0}^L\rho |v_t|^2{dx}+\left(\sqrt{\frac{\rho}{\alpha_1}}+\sqrt{\frac{\mu\gamma^2}{\alpha_1}}\right)\int_{0}^L\alpha_1 |v_x|^2{dx}\right.\\
			&  +\left(\sqrt{\frac{\mu}{\beta}}+\sqrt{\frac{\mu\gamma^2}{\alpha_1}}\right)\int_{0}^L\mu |p_t|^2{dx}\left.+\sqrt{\frac{\mu}{\beta}}\int_{0}^L\beta |\gamma v_x-p_x|^2{dx}\right]\\
			&\leq L\eta E(t).
		\end{array}
	\end{equation*}
	Since { $|F(t)|\leq L\eta E(t),$} this leads to
	$\left|E_\delta(t)\right|  \leq \left|E(t)\right|+\delta \left|F(t) \right| \le \left(1+L\eta \delta  \right) E(t),$
	and analogously, $\left|E_\delta(t)\right|  \ge \left(1-L\eta \delta  \right) E(t), $ and therefore, \eqref{equiva} is immediate
\end{proof}
\begin{Lem}\label{lem1}
	For any $\e, t>0$, $F(t)$ in \eqref{eq5} satisfies the following inequalities,
   \begin{eqnarray}\label{eqx}
	\begin{array}{ll}	
        &\frac{dF}{dt}\leq-E(t)+\frac{L}{2}\left[\rho+\frac{(1+\e)\xi_1^2}{\alpha_1}\right] \left|v_t(L,t)\right|^2\\
		&\qquad\quad+\frac{L}{2}\left[\mu+{\left(1+\frac{1}{\e}\right)}\frac{\xi_2^2\gamma^2}{\alpha_1}+\frac{\xi_2^2}{\beta}\right]\left|p_t(L,t)\right|^2.
    \end{array}
   \end{eqnarray}
\end{Lem}
\noindent\begin{proof}
	Recalling $\alpha_1=\alpha-\gamma^2\beta$, and \eqref{eq1},
	\begin{eqnarray*}
		\begin{array}{l}
			\frac{dF}{dt}		=\frac{L}{2}\beta\left(\gamma v_x(L,t)-p_x(L,t)\right)^2+\frac{L}{2}\alpha_1 (v_x(L,t))^2\\
		\qquad+\frac{L}{2}\left[\rho \left|v_t(L,t)\right|^2+\mu \left|p_t(L,t)\right|^2\right]-E(t).
		\end{array}
	\end{eqnarray*}
	Next, the boundary conditions are used so that
\begin{eqnarray*}
\begin{array}{l}
	\frac{dF}{dt}=-E(t)+\frac{L}{2}\left(\mu+\frac{\xi_2^2}{\beta}\right)\left|p_t(L,t)\right|^2\\
	\quad+\frac{L\rho}{2}\left|v_t(L,t)\right|^2+\frac{L}{2\alpha_1}(\xi_1\left|v_t(L,t)\right|+\xi_2\gamma{\left|p_t(L,t)\right|)^2}.
\end{array}
\end{eqnarray*}
	Finally, by the {generalized} Young's inequality { with $\e$} (or  Peter–Paul inequality),  \eqref{eqx} is obtained for any $\e>0$.
\end{proof}
Now, the exponential stability result takes the following form
\begin{thm}\label{thm2}
	The energy $E(t)$ of solutions decays exponentially, \ie {for any $\e >0$,}
	\begin{equation}\label{eq:thm2}
		\begin{array}{ll}
			&E(t)\leq M E(0)e^{-\sigma t},\qquad \forall t>0\\
			&\sigma(\delta)= \delta\left(1-\delta L \eta\right), \qquad M(\delta)=\frac{1+\delta L \eta}{1-\delta L\eta},
		\end{array}
	\end{equation}
    \begin{equation}\label{casei}
        \begin{cases}
            &{\delta(\xi_1,\xi_2,\epsilon)<\frac{1}{L} \min\left(\frac{1}{\eta}, f_1(\xi_1,\epsilon),f_2(\xi_2,\epsilon)\right),}\\
			&f_1(\xi_1,\epsilon):=\frac{2\xi_1\alpha_1}{\rho\alpha_1+(1+\epsilon)\xi_1^2},\\
			&f_2(\xi_2,\epsilon):=\frac{2\xi_2\e\alpha_1\beta}{\e\mu\alpha_1\beta+(\e\alpha+\gamma^2\beta)\xi_2^2}.
        \end{cases}
    \end{equation}
\end{thm}
\begin{proof} Considering $\frac{d E_\delta}{dt}=\frac{dE}{dt}+\delta\frac{dF}{dt}$ with Lemma \ref{lem1}
	\begin{align*}
		&\frac{d E_\delta}{dt}\leq \underbrace{\left[\frac{\delta L}{2}\left[\rho+\frac{(1+\e)\xi_1^2}{\alpha_1}\right]-\xi_1\right]}_{\le 0}\left|v_t(L,t)\right|^2\\
		&+\underbrace{\left[\frac{\delta  L}{2}\left[\mu+\frac{(\e\alpha +\gamma^2\beta)\xi_2^2}{\e\alpha_1\beta }\right]-\xi_2\right]}_{\le 0}\left|p_t(L,t)\right|^2-\delta E(t)
	\end{align*}
	where $\delta$ is chosen to make the coefficients  nonpositive, \ie
	$
		\delta <\frac{1}{L} \min\left(\frac{1}{\eta}, f_1(\xi_1,\epsilon),f_2(\xi_2,\epsilon)\right).
    $
	Next, by the equivalence of $E(t)$ and $E_\delta(t)$ from Lemma \ref{lem2},
	\begin{eqnarray}\label{21}
		&&\frac{dE_\delta}{dt}\leq- \delta\left(1-\delta L \eta\right)E_\delta(t).
	\end{eqnarray}
	By choosing $\sigma=\delta\left(1-\delta L \eta\right)>0,$
	together with \eqref{21} lead to
	\begin{equation}\label{exp}
		E_\delta (t) \leq E_\delta (0) e^{-\sigma t}.
	\end{equation}
	Hence, \eqref{exp} with Lemma \ref{lem2} lead to the desired result.
\end{proof}

\section{Optimization of Sensor Feedback Amplifiers}
\label{sec:OptimSensor}
{ The decay rate $-\sigma$ in Theorem~\ref{thm2} provides an upper
 bound for the exponential decay rate of the system. In this section, we provide an optimization process for feedback amplifiers $\xi_1$ and $\xi_2$ that ensure the maximal value of $\sigma$, given by
\begin{equation}\label{maximal}
	\sigma_{max}(\delta)=\frac{1}{4\eta L} ~{\rm achieved ~at}~ \delta=\frac{1}{2\eta L}.
\end{equation}
}
Since $\delta$ and $\sigma$ are functions of $\epsilon$ and the sensor feedback amplifiers $\xi_1$ and $\xi_2$ in~\eqref{c1}, the following results provide intervals for the feedback amplifiers and $\epsilon$ that ensures the decay rate~\eqref{maximal} is attained.

\begin{thm}\label{theorem4} Define non-negative constants
	\begin{equation}
		\label{cs}\begin{array}{ll}
			c_1^\pm:=\frac{{ 2}\alpha_1\eta\pm\sqrt{{ 4}\alpha_1^2\eta^2-(1+\epsilon)\rho\alpha_1}}{1+\epsilon},\\
			c_2^\pm:=\frac{{ 2}\epsilon\alpha_1\beta\eta\pm\sqrt{{ 4}(\epsilon\alpha_1\beta\eta)^2-(\epsilon\alpha+\beta\gamma^2)\epsilon\mu\alpha_1\beta}}{\epsilon\alpha+\beta\gamma^2},
		\end{array}
	\end{equation}
	with any $\e$ such that
	\begin{equation}\label{epsilon_bound1}
		\frac{\beta\gamma^2\mu}{{ 4}\alpha_1\beta\eta^2-\alpha\mu}{ <}\epsilon{ <}\frac{{ 4}\alpha_1\eta^2-\rho}{\rho}.
	\end{equation}
	
	The  decay rate $-\sigma_{max}(\delta)$ is achieved for the closed-loop system \eqref{eq1},\eqref{c1} when the feedback amplifiers are chosen such as
	$\xi_1\in { (}c_1^-, c_1^+{ )},  \xi_2\in { (}c_2^-,c_2^+{ )}.$
\end{thm}
\begin{proof}
	Observe that to achieve the  $\sigma_{max}(\delta)$ in \eqref{maximal}, it is sufficient to have $f_1(\xi_1,\epsilon){ >}\frac{1}{{ 2}\eta}$ and $f_2(\xi_2,\epsilon){ >}\frac{1}{{ 2}\eta}$.
	
	Thus, $f_1(\xi_1,\epsilon){ >} \frac{1}{{ 2}\eta},$ \ie $\frac{2\xi_1\alpha_1}{\rho\alpha_1+(1+\epsilon)\xi_1^2}{ >} \frac{1}{{ 2}\eta},$
	implies that
	\begin{eqnarray}\label{upperbound}\begin{array}{ll}
		\e { <} h_1(\xi_1):=\frac{{ 4}\alpha_1\xi_1\eta-\rho\alpha_1-\xi_1^2}{\xi_1^2}.
	\end{array}
	\end{eqnarray}
	
	Noting that ${ 4}\alpha_1^2\eta^2-\rho\alpha_1>0$ by \eqref{eta}, $h_1(\xi_1)$ defines an upper bound for $\e$, and observe that $\xi_1$ must be chosen in between the following roots of $h_1$ to ensure $\e>0$, see Fig. \ref{fig1}
	\begin{eqnarray*}
		\begin{array}{l}
			a_1^\pm:={ 2}\alpha_1\eta\pm\sqrt{{ 4}\alpha_1^2\eta^2-\rho\alpha_1},
		\end{array}
	\end{eqnarray*}
	Observe that $h_1(\xi_1)\geq0$ if and only if $\xi_1\in(a_1^-,a_1^+)$. Seeking the critical points of $h_1(\xi_1),$
	$\frac{\partial h_1}{\partial \xi_1}=\frac{-{ 4}\xi_1^2\alpha_1\eta+2\xi_1\rho \alpha_1}{\xi^4}=0,$
	leads to $\xi_1=\frac{\rho}{{ 2}\eta},$ for which $h_1$ achieves its maximum value. Substituting $h_1\left(\xi_1=\frac{\rho}{\eta} \right)$ into \eqref{upperbound} yields the following upper bound for $\e$
	\begin{equation}\label{epsilonupper}
		\e { <} \frac{{ 4}\alpha_1\eta^2-\rho}{\rho}.
	\end{equation}
	By $f_2(\xi_2,\epsilon){ >} \frac{1}{{ 2}\eta},$
$
			\epsilon{ >} h_2(\xi_2):=\frac{\beta\gamma^2\xi_2^2}{{ 4}\alpha_1\beta\xi_2\eta-\mu\alpha_1\beta-\alpha\xi_2^2}.
		$	Since $\e >0$, the denominator ${ 4}\alpha_1\beta\xi_2\eta-\mu\alpha_1\beta-\alpha\xi_2^2$ is chosen to be strictly positive. This leads to $\xi_2\in (a_2^-,a_2^+)$ where
	\begin{eqnarray}	
		\begin{array}{ll}
			a_2^\pm:=&\frac{{ 2}\alpha_1\beta \eta\pm\sqrt{{ 4}\alpha_1^2\beta^2\eta^2-\alpha\mu\alpha_1\beta}}{\alpha},
		\end{array}
	\end{eqnarray}
	and $\xi_2=a_2^-$ and $\xi_2=a_2^+$ are the two vertical asymptotes of $h_2(\xi_2),$ see dashed lines in Fig. \ref{fig1}. Note that this condition ensures $h_2(\xi_2)\geq0$. Seeking the critical points of $h_2(\xi_2)$
	leads to  $\xi_2=\frac{\mu}{\eta},$ for which $h_2(\xi_2)$ takes its minimum value. Substituting $h_2\left(\frac{\mu}{{ 2}\eta} \right)$ into \eqref{upperbound} yields the lower bound for $\e$
	\begin{eqnarray}\label{epsilonlower}
		\begin{array}{l}
			\e { >}\frac{\beta\gamma^2\mu}{{ 4}\alpha_1\beta\eta^2-\alpha\mu}
		\end{array}
	\end{eqnarray}
	where ${ 4}\alpha_1\beta\eta^2-\alpha\mu>0$ by \eqref{eta}. Restricting $h_1(\xi_1)$ and $h_2(\xi_2)$ in between roots and asymptotes, respectively, $\e$ values corresponding to the filled regions in Fig. \ref{fig1} satisfy conditions \eqref{epsilonupper} and \eqref{epsilonlower}, respectively.
	
	\noindent For any $\e$ that satisfies \eqref{epsilon_bound1}, $f_1(\xi_1){ >}\frac{1}{{ 2}\eta}$, and $f_2{ >}\frac{1}{{ 2}\eta}$, \ie
	\begin{equation*}
		\begin{array}{ll}
			-(1+\e)\xi_1^2+{ 4}\alpha_1\eta\xi_1-\rho\alpha_1\geq0,\\
			-(\e\alpha+\beta\gamma^2)\xi_2^2+{ 4}\e\alpha_1\beta\eta\xi_2-\e\mu\alpha_1\beta\geq0,
		\end{array}
	\end{equation*}
	yield $\xi_1\in { (}c_1^-,c_1^+{ )}$ and $\xi_2\in { (}c_2^-,c_2^+{ )}$, respectively.
	
\begin{figure}[h!]
	\centering
	\includegraphics[width = 0.22\textwidth ]{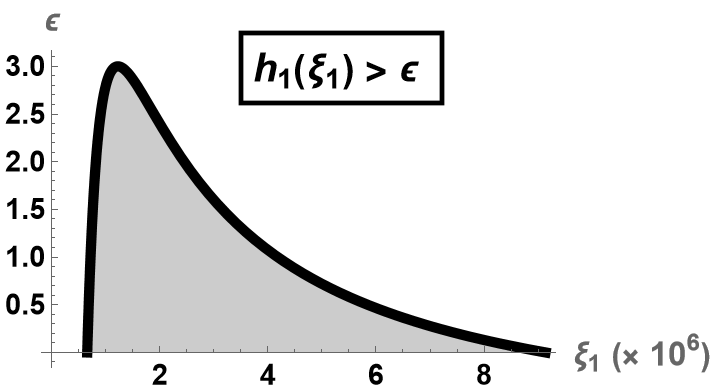} \includegraphics[width = 0.22\textwidth ]{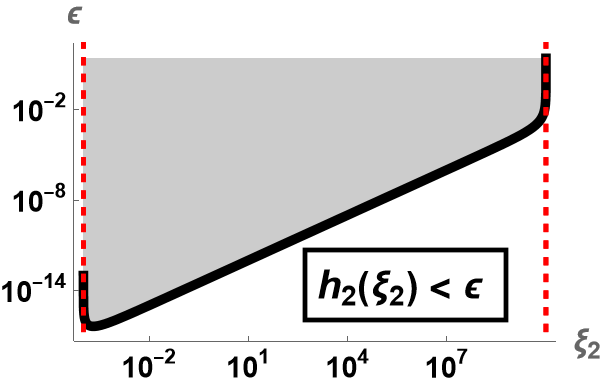}
		\caption{For the material parameters in Table I, bounds of $\e$ depend on $h_1(\xi_1)$ and $h_2(\xi_2)$. Note that the $h_2-$plot is logarithmic.}
	\label{fig1}
\end{figure}

	We prove that such $\xi_1\in { (}c_1^-,c_1^+{ )}$ and $\xi_2\in { (}c_2^-,c_2^+{ )}$ exist. Note that the sufficient condition for ${ (}c_1^-,c_1^+{ )}, { (}c_2^-,c_2^+{ )}\neq \emptyset$ is the existence of $\e$ satisfying \eqref{epsilon_bound1} which is ensured by
	\begin{eqnarray}\label{epsilon_bound}
		\begin{array}{ll}
			\frac{4\alpha_1\eta^2-\rho}{\rho}\geq \frac{3\alpha_1\eta^2}{\rho} \geq 3, &\frac{\beta\gamma^2\mu}{{ 4}\alpha_1\beta\eta^2-\alpha\mu} \leq \frac{\beta\gamma^2\mu}{3\alpha\mu} \leq \frac{1}{3}.
		\end{array}
	\end{eqnarray}
\end{proof}
\begin{figure}[h!]
	\centering
	\includegraphics[width=0.46\columnwidth ]{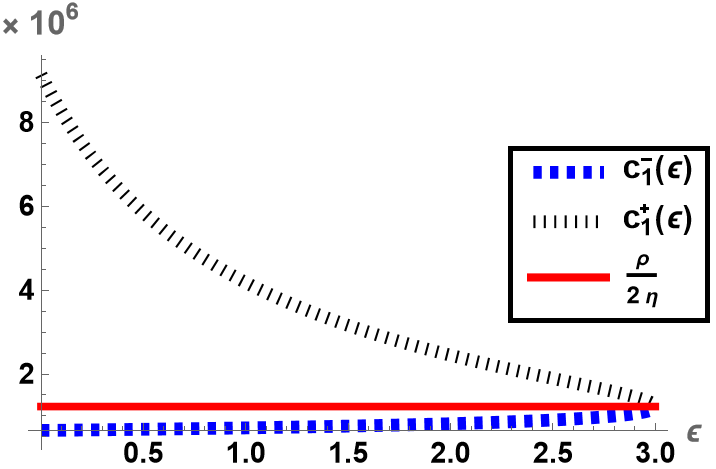} \includegraphics[width=0.52\columnwidth ]{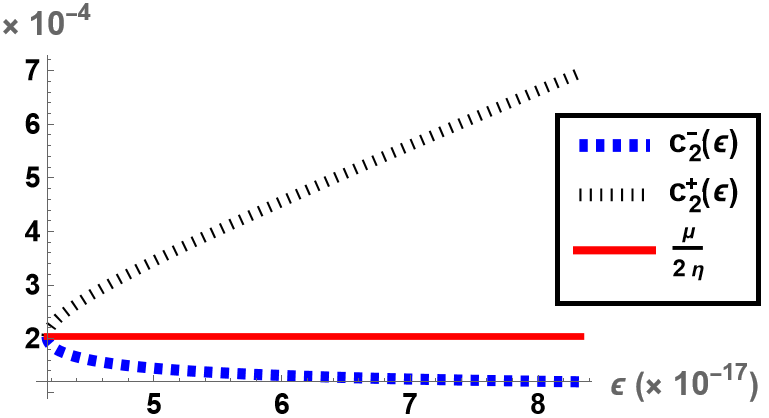}
	\caption{For the material parameters  in Table I, graphs of $c_1^\pm,c_2^\pm$ with respect to $\e$ satisfy the inequalities in \eqref{epsilon_bound1}.  }
	\label{figb11}
\end{figure}
The amplifier $\xi_1$  maximizing $h_1$ and the amplifier $\xi_2$ minimizing $h_2$ are always in the respective intervals, \ie $\xi_1=\frac{\rho}{{ 2}\eta}\in{ (}c_1^-,c_1^+{ )}$ and $\xi_2=\frac{\mu}{{ 2}\eta}\in { (}c_2^-,c_2^+{ )}$. Moreover, as $\e$ approaches to the upper (lower) bound ${ (}c_1^-,c_1^+{ )}$ gets smaller (larger) and ${ (}c_2^-,c_2^+{ )}$ gets larger (smaller), see Fig. \ref{figb11}. Feedback amplifiers chosen within the intervals ensure that $f_1(\xi_1)\geq\frac{1}{{ 2}\eta}$ and $f_2(\xi_2)\geq\frac{1}{{ 2}\eta}$ are satisfied, see Figs.~\ref{figb12a}~and~\ref{figb12b}.

\begin{figure}[h!]
	\centering
	\includegraphics[width=0.8\columnwidth ]{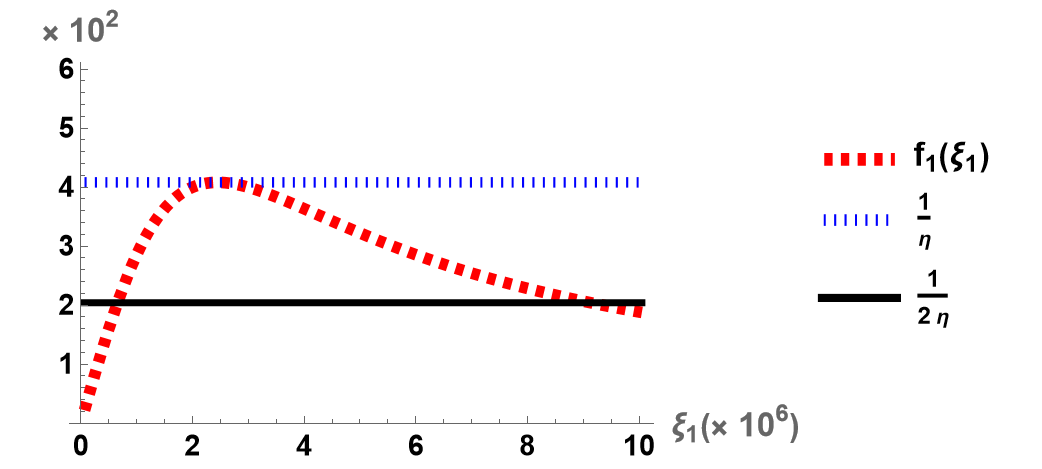}
	\caption{For $\epsilon\approx 10^{-15},$ $f_1(\xi_1)>\frac{1}{{ 2}\eta}$ when  $\xi_1\in { (}c_1^-, c_1^+{ )}.$}
	\label{figb12a}
\end{figure}
\begin{figure}[h!]
	\centering
\includegraphics[width=0.8\columnwidth ]{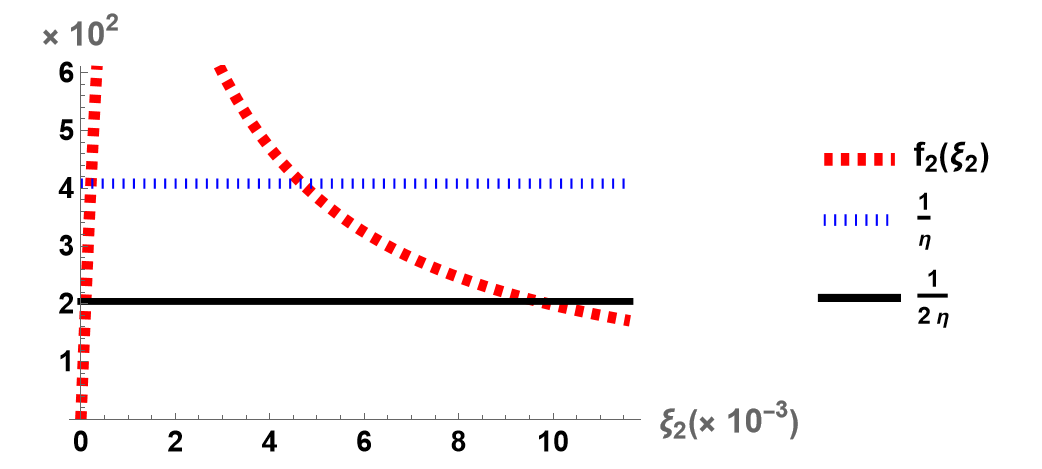}
	\caption{For $\epsilon\approx 10^{-15},$ $f_2(\xi_2)>\frac{1}{{ 2}\eta}$ when  $\xi_2\in { (}c_2^-,c_2^+{ )}.$}
	\label{figb12b}
\end{figure}

\section{Simulations and Numerical Experiments}
{
In this section, we present numerical simulations primarily aimed at demonstrating the following:
\begin{itemize}
	\item The exponential stability of the system \eqref{eq1} with \eqref{c1} is showcased using  feedback amplifiers $\xi_1\in (c_1^-, c_1^+)$ and $\xi_2\in(c_2^-,c_2^+)$, as obtained in Theorem~\ref{theorem4}.
	\item Comparison of the reduced model exponential decay rates for a selection of feedback amplifiers.
\end{itemize}	

Moreover, we address a common challenge encountered in standard model reductions, such as Finite Differences (FD) and Finite Elements (FE), which introduce spurious high-frequency vibrational modes to the system. These spurious modes significantly impact the decay rate of the system. In fact, reduced models lack the uniform observability property as the mesh parameter approaches zero, necessitating filtering techniques such as the direct Fourier filtering~\cite{B, I-Z} or the indirect filtering \cite{T-Z}.

To mitigate the adverse effects of approximation methods on the decay rate, we consider a recently introduced technique known as Order-Reduced Finite Differences (ORFD)~\cite{GuoWave}. The ORFD method has been observed to retain the decay rate of the infinite-dimensional system. For rigorous stability analyses of all three approximations, interested readers can refer to \cite{O-Emran2}.}

Consider a given natural number $N\in\mathbb{N}$, which denotes the number of nodes in the spatial semi-discretization. Let's introduce the mesh size as $h:=\frac{1}{N+1}$, and discretize the interval $[0,L]$ as
$0=x_0<x_1<\dots <x_j=j*h<\ldots<x_N<x_{N+1}=L.$
Let $(v_j,p_j)=(v_j,p_j)(t)\approx (v,p)(x_j,t)$ be the approximation of the solution $(v,p)(x,t)$ of \eqref{eq1}-\eqref{c1} at the point space $x_j=j\cdot h$ for any $j=0,1,...,N,N+1$, and let $\vec{v}=[v_1,v_2,...,{v_{N+1}}]^T$ and  $\vec{p}=[p_1,p_2,...,{p_{N+1}}]^T$.

{In the ORFD method, in addition to the standard nodes $\left\{x_{j}\right\}_{j=0}^{N},$ we also consider the in-between middle nodes of each subinterval, denoted by $\left\{x_{j+\frac{1}{2}}:=\frac{x_{j+1}+x_j}{2}\right\}_{j=0}^{N}.$ Define the average, $v_{j+\frac{1}{2}}:=\frac{v_{j+1}+u_j}{2}$ and difference operators
\begin{eqnarray}
\begin{array}{ll}
 \delta_x v_{j+\frac{1}{2}}:=\frac{v_{j+1}-v_j}{h},&\delta_x^2 v_{j}:=\frac{v_{j+1}-2v_j+v_{j-1}}{h^2},
\end{array}
\end{eqnarray}
It is worth noting that by considering odd-number derivatives at the in-between nodes within the uniform discretization of $[0,L]$, higher-order approximations are achieved \cite{GuoWave}.

The ORFD approximation of equations \eqref{eq1} with \eqref{c1} is
\begin{equation}\label{ORFD}
		\left(\bm C_1 \otimes \bm M\right) \begin{bmatrix}
			{\ddot{\vec{v}}}   \\
			{\ddot{\vec{p}}}        \\
		\end{bmatrix} +\left(\bm C_2  \otimes \bm A_h\right) \begin{bmatrix}
			{{\vec{v}}}   \\
			{{\vec{p}}}        \\
		\end{bmatrix} \\
		 + \left(\bm C_3 \otimes  \mathcal{B}\right) \begin{bmatrix}
			{\dot{\vec{v}}}   \\
			{\dot{\vec{p}}}        \\
		\end{bmatrix}
		=\vec 0,
\end{equation}
where $C_1$ and $C_2$ are matrices for material parameters, and $C_3$ is the matrix for the state feedback amplifiers, defined as
\[
	\bm C_1 = \begin{bmatrix}
		\rho & 0 \\
		0 & \mu
	\end{bmatrix},\quad \bm C_2 = \begin{bmatrix}
		\alpha & -\gamma\beta \\
		-\gamma\beta & \beta
	\end{bmatrix}, \quad \bm C_3 = \begin{bmatrix}
		\xi_1 & 0 \\
		0 & \xi_2
	\end{bmatrix},
\]
the $(N+1)\times (N+1)$ mass matrix $\bm M$ is defined by
\[
	\bm M  =\frac{1}{4}\begin{bmatrix}
		2&1&0&\dots&\dots&\dots&0\\
		1&2&1&0&\dots&\dots&0\\
		&\ddots&\ddots&\ddots&\ddots&\ddots&\\
		0&\dots&\dots&0&1&2&1\\
		0&\dots&\dots&\dots&0&1&1\\
	\end{bmatrix},
\]
the $(N+1)\times (N+1)$ central difference matrix $\bm A_h$ is
\[
	\bm A_h=\frac{1}{h^2}\begin{bmatrix}
		2&-1&0&\dots&\dots&\dots&0\\
		-1&2&-1&0&\dots&\dots&0\\
		&\ddots&\ddots&\ddots&\ddots&\ddots&\\
		0&\dots&\dots&0&-1&2&-1\\
		0&\dots&\dots&\dots&0&-1&1\\
	\end{bmatrix},
\]
and the $(N+1)\times (N+1)$ boundary node matrix $\mathcal{B}$ is a zero matrix except for the $(N+1)\times (N+1)$-th entry, which is $\frac{1}{h},$ due to the boundary damping being injected at the last node. Here, $\otimes$ denotes the matrix Kronecker product. It's worth noting that equation \eqref{ORFD} can be rewritten in the first-order form as
\begin{equation}
	\label{ORFDFirstOrder}
	\frac{d}{dt}\begin{bmatrix}
		{{\vec{v}}}   &
		{{\vec{p}}}    &
		{\dot{\vec{v}}}   &
		{\dot{\vec{p}}}
	\end{bmatrix}^T =
	 \bm{\mathcal{A}}\begin{bmatrix}
		{{\vec{v}}}   &
		{{\vec{p}}}    &
		{\dot{\vec{v}}}   &
		{\dot{\vec{p}}}
	\end{bmatrix}^T,
\end{equation}
\begin{equation}\label{calA}
	\bm{\mathcal{A}}=\begin{bmatrix}
		\bm 0_{2N+2} & I_{2N+2}\\
				-\left(\bm C_1^{-1}\bm C_2\right) \otimes \left(\bm M^{-1}\bm A_h\right) & -\left(\bm C_1^{-1}\bm C_3\right) \otimes \left(\bm M^{-1}\mathcal{B}\right)
			\end{bmatrix}
\end{equation}
The discretized  energy corresponding to \eqref{ORFDFirstOrder} is
\begin{equation}\label{ORFFDEnergy}
	\begin{split}
		E_h(t):=&\frac{h}{2}\left\langle (\bm C_1 \otimes \bm M)  \begin{bmatrix}
			{\dot {\vec v}} \\
			{\dot {\vec p}}
		\end{bmatrix},\begin{bmatrix}
			{\dot {\vec v}} \\
			{\dot {\vec p}}
		\end{bmatrix}\right\rangle\\
		&\qquad\qquad +\frac{h}{2}\left\langle (\bm C_1 \otimes \bm A_h)  \begin{bmatrix}
			{ {\vec v}} \\
			{ {\vec p}}
		\end{bmatrix},\begin{bmatrix}
			{ {\vec v}} \\
			{ {\vec p}}
		\end{bmatrix}\right\rangle.
	\end{split}
\end{equation}

	To show the importance of our analysis for the choices of sensor feedback amplifiers via the optimization process outlined in the previous section, sample numerical simulations are presented for the material constants in Table \ref{table}. }
		\begin{table}[htb!]
		\centering
		\caption{Realistic Piezoelectric Material Parameters}
		\label{table}
		\begin{tabular}{c|c|c|c}
			\hline
			Parameter&Symbol&Value&Unit\\
			\hline
			Length of the beam&$L$&$1$&m \\
			Mass density&$\rho$&$6000$&kg/${\rm m}^3$\\
			Magnetic permeability&$\mu$&$10^{-6}$ & H/m\\
			Elastic stiffness&$\alpha$&$10^9$& N/${\rm m}^2$\\
			Piezoelectric constant&	$\gamma$&$10^{-3}$&C/${\rm m}^3$\\
			Impermittivity  &$\beta$& $ 10^{12}$& m/F\\
			\hline
		\end{tabular}
	\end{table}
	
	{ Note that $\sigma_{max}(\delta)$ in  \eqref{maximal} depends solely on the material parameters. In particular, $\sigma_{\rm max} =\frac{1}{4\eta L}\approx 102.04$ for the material parameters in Table \ref{table}. Letting $\epsilon = 1$, which is a valid choice for any system according to equations \eqref{epsilon_bound}, the intervals for optimal sensor feedback amplifiers in \eqref{cs}, corresponding to  $\sigma_{\rm max},$ are as follows:
		\begin{equation}\label{optimal_feedback_amplifiers}
		\begin{split}
			(c_1^-,c_1^+) &\approx (7.17\times 10^5,    ~4.18\times 10^6),\\
			(c_2^-,c_2^+) &\approx (1.02\times 10^{-4}, ~9.78\times 10^{9}).
		\end{split}
	\end{equation}
Here note that  both the maximal decay rate  $-\sigma_{\rm max}$ and the corresponding intervals  for feedback amplifiers in \eqref{optimal_feedback_amplifiers} are independent of the choice of initial conditions.

For simulations,  we take $N=80$ nodes, and we choose  triangular  (hat-type) initial conditions that are composed of high and low-frequency vibrational nodes. The simulations are performed for $T=0.1$ sec to show the rapid exponential decay of solutions. Indeed, the solutions of the system \eqref{ORFD} with the choices of optimal feedback amplifiers $\xi_1 =10^6 \in (c_1^-,c_1^+)$ and $\xi_2 =10^9 \in (c_2^-, c_2^+)$, are observed drive all states to zero state rapidly, e.g. see Fig.~\ref{optimal3Dplot}, and the total energy follows the exponentially decay in a higher rate than claimed $\sigma_{\rm max}$, Fig.~\ref{energies}.

Our analysis indicates that even a relatively small perturbation in the values of the feedback amplifiers $\xi_1$ can lead to suboptimal performance. Indeed, when $\xi_1 = 10^4$ and does not fall within the range $(c_1^-, c_1^+)$, the system still continues to exhibit low and high-frequency longitudinal oscillations for $T=0.1$ seconds, as demonstrated in Figure~\ref{NonOptimal3Dplot}. Furthermore, the energy of the system decays in an unstable manner, with a decay rate significantly larger than $-\sigma_{\rm max}$, as depicted in Figure~\ref{energies}. }

\begin{figure}[htb!]
	\centering
	\includegraphics[width=0.8\columnwidth]{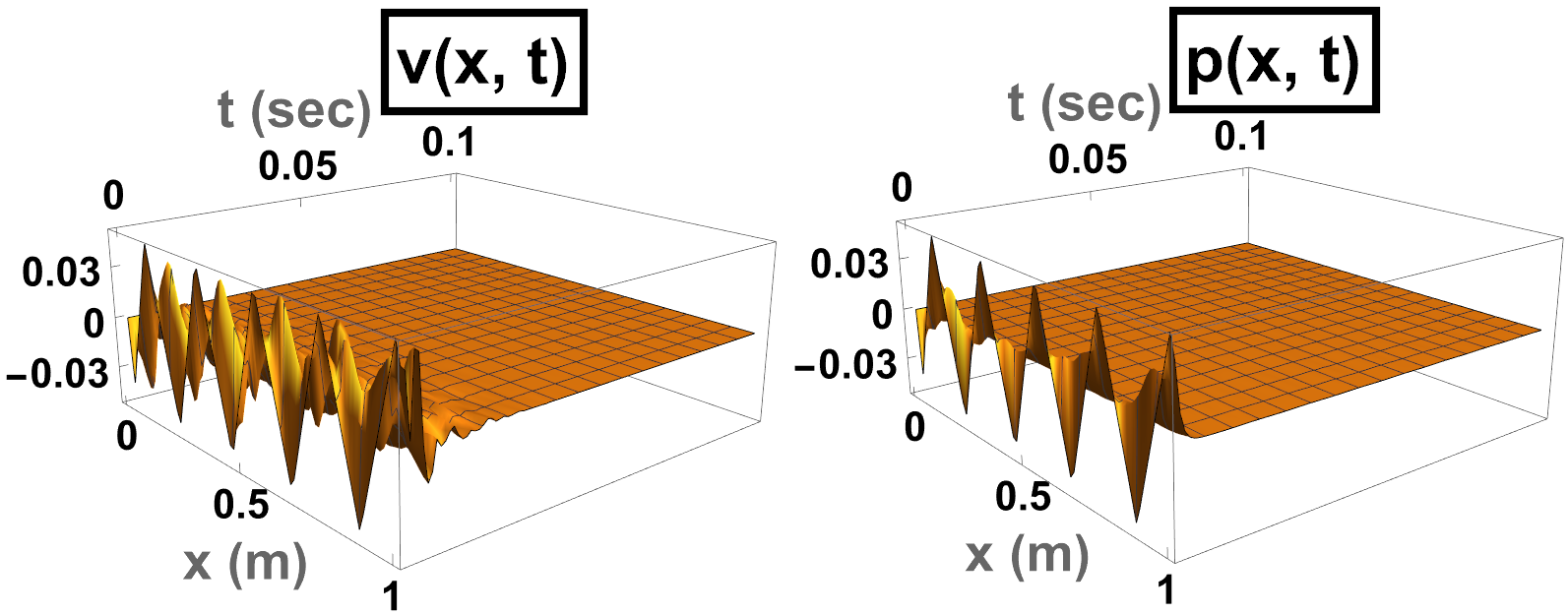}
	\caption{System behavior for the optimal choice of feedback amplifiers: $\xi_1 = 10^6$ and $\xi_2 = 10^{9}$. Both solutions $v(x,t)$ and $p(x,t)$ of Equation~\eqref{ORFD} decay rapidly to the zero state in much less than $0.1$ sec.}
	\label{optimal3Dplot}
\end{figure}

\begin{figure}[htb!]
	\centering
	\includegraphics[width=0.8\columnwidth]{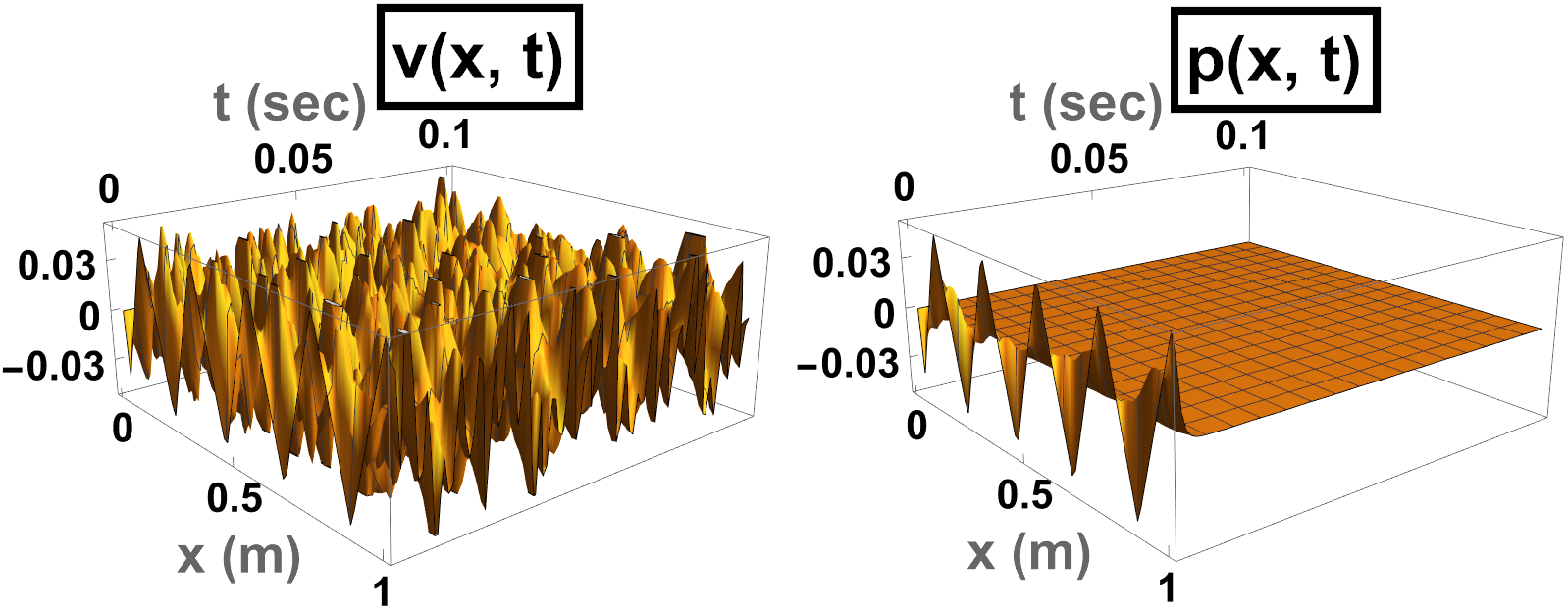}
	\caption{For a non-optimal choice of $\xi_1 = 10^4$, which control longitudinal strains, and an optimal choice of $\xi_2 = 10^{9}$, controlling the voltage, the longitudinal oscillations $v(x, t)$ for eq.~\eqref{ORFD} exhibit a significantly slower decay rate.}
	\label{NonOptimal3Dplot}
\end{figure}

    \begin{table}
        \centering
        \begin{tabular}{c|c c c c c c}
            \backslashbox{$\xi_1$}{$\xi_2$}& $10^{5}$ & $10^{5.5}$ &\tikzmark{top left 1}		$10^{6}$ & $10^{6.5}$&  $10^{7}$&  $10^{7.5}$\\
            \hline
            $ 10^{-7}$  & $-0.1$ & $-0.1$ &$ -0.1 $& $-0.1$& $-0.1$&$-0.1$\\
            $ 10^{-5}$ & $-10$ & $-10$ &$ -10$& $-10$& $-10$& $-10$\\
            \tikzmark{top left 2} $ 10^{-3}$ & $-17$ & $-53$ &$-177 $& $-421$& $-101.9$&$-32$\\
            $ 10^{6}$& $-17$ & $-53$ &$-177 $& $-421 $& $-101.9$&$-31$\\
            $ 10^{9}$& $-17$ & $-53$ &$-177 $& $-421 $& $-101$&$-30$\tikzmark{bottom right 2}\\
            $ 10^{10}$& $-17$ & $-52$ &$-100 $& $-100 $& $-97$&$-23$\\
            $ 10^{11}$ & $-10$ & $-10$ &$-10 $& $-10$\tikzmark{bottom right 1}& $-10$&$-9$
        \end{tabular}
        \caption{The $\max\{\Re(\mu_k)\}$ values for different choices of feedback amplifiers $\xi_1$ and $\xi_2$. Red boxes indicate the intervals $(c_1^-, c_1^+)$ and $(c_2^-,c_2^+)$, where the maximal decay rate is achieved.}
        \label{DecayTable}
        \DrawBox[ultra thick, red]{top left 1}{bottom right 1}
        \DrawBox[ultra thick, dashed, red]{top left 2}{bottom right 2}
    \end{table}

\noindent\indent{
To demonstrate the robustness of our theoretical findings, we examine the spectrum of the numerical approximation. Let $\mu_k$ denote the eigenvalues of $\bm{\mathcal{A}}$ as defined in Equation \eqref{calA}. The decay rate of the system \eqref{eq1} with \eqref{c1}, and the theoretically derived decay rate $-\sigma$ in Theorem \eqref{thm2}, can be approximated as the maximal real part of the eigenvalues of $\bm{\mathcal{A}}$, specifically $-\sigma \approx \max\{\Re(\mu_k)\}$.
A contour plot of $\max\{\Re(\mu_k)\}$ in terms of the feedback amplifiers $\xi_1$ and $\xi_2$ is presented in Figure \ref{DecayRates}. It is notable that as the feedback amplifiers $\xi_1$ and $\xi_2$ fall within the intervals $(c_1^-, c_1^+)$ and $(c_2^-,c_2^+)$, respectively shown as red lines, the maximal decay rate is achieved. Conversely, the decay rate significantly diminishes outside of these intervals. Exact  values of $\max\{\Re(\mu_k)\}$ for selected feedback amplifiers are also provided in Table \ref{DecayTable}. As a final note, the decay rates are much smaller than $-\sigma_{\rm max}$ when both feedback amplifiers are chosen from the intervals $(c_1^-, c_1^+)$ and $(c_2^-,c_2^+)$.
}

\color{black}
	\section{Conclusions \& Future Work}

In summary, using the Lyapunov approach, we derived the exponential decay rate $``-\sigma"$ for system \eqref{eq1}, with the design of state feedback amplifiers $\xi_1$ and $\xi_2$ as described in \eqref{c1}. This decay rate ensures fast stabilization, though improvements are possible. The rate depends on material parameters and feedback amplifiers. Our numerical results confirm the robustness of this decay rate, with deviations affecting it from the theoretical value. Interactive simulations are available through the Wolfram Demonstrations Project~\cite{W}.

A full numerical analysis and proof of exponential stability for the semi-discretized model \eqref{ORFD}, as $h \to 0$, and determining optimal designs of feedback amplifiers remain beyond this paper's scope and are ongoing research. However, the Lyapunov approach here serves as the basis for results in~\cite{O-Emran2}, where the discretized model achieves the same decay rate as the PDE model.

Future work can extend these concepts to multi-layer magnetizable piezoelectric beams, which have numerous practical applications \cite{O3}. The design of feedback amplifiers may become more delicate when there are more than two state feedback amplifiers, as highlighted in \cite{O3}. There is also potential to explore general hyperbolic PDE systems with multiple boundary dampers. Relevant studies address systems with random inputs \cite{Marin2017} and the optimal design and placement of actuators in higher-dimensional systems \cite{PeriagoMunch,Mor-Dem}, both of which are valuable for future research.

\begin{figure}[htb!]
	\centering
	\includegraphics[width=0.9\columnwidth]{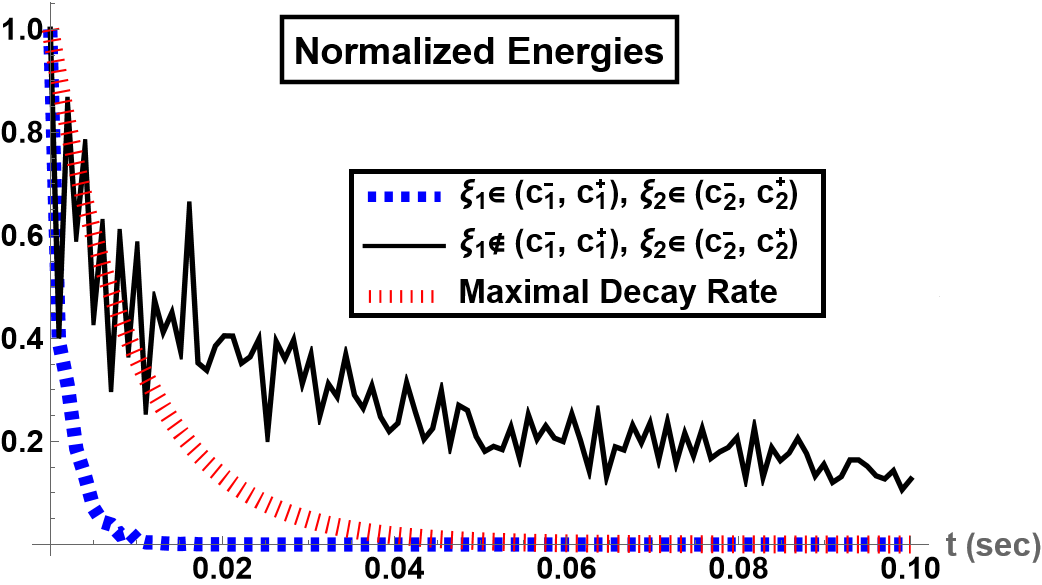}
	\caption{Normalized total energies for \eqref{ORFD} with optimal and non-optimal feedback amplifiers, along with the theoretically derived energy decay as stated in Theorem~\eqref{theorem4}.}
	\label{energies}
\end{figure}	

\begin{figure}[htb!]
	\centering
	\includegraphics[width=\columnwidth]{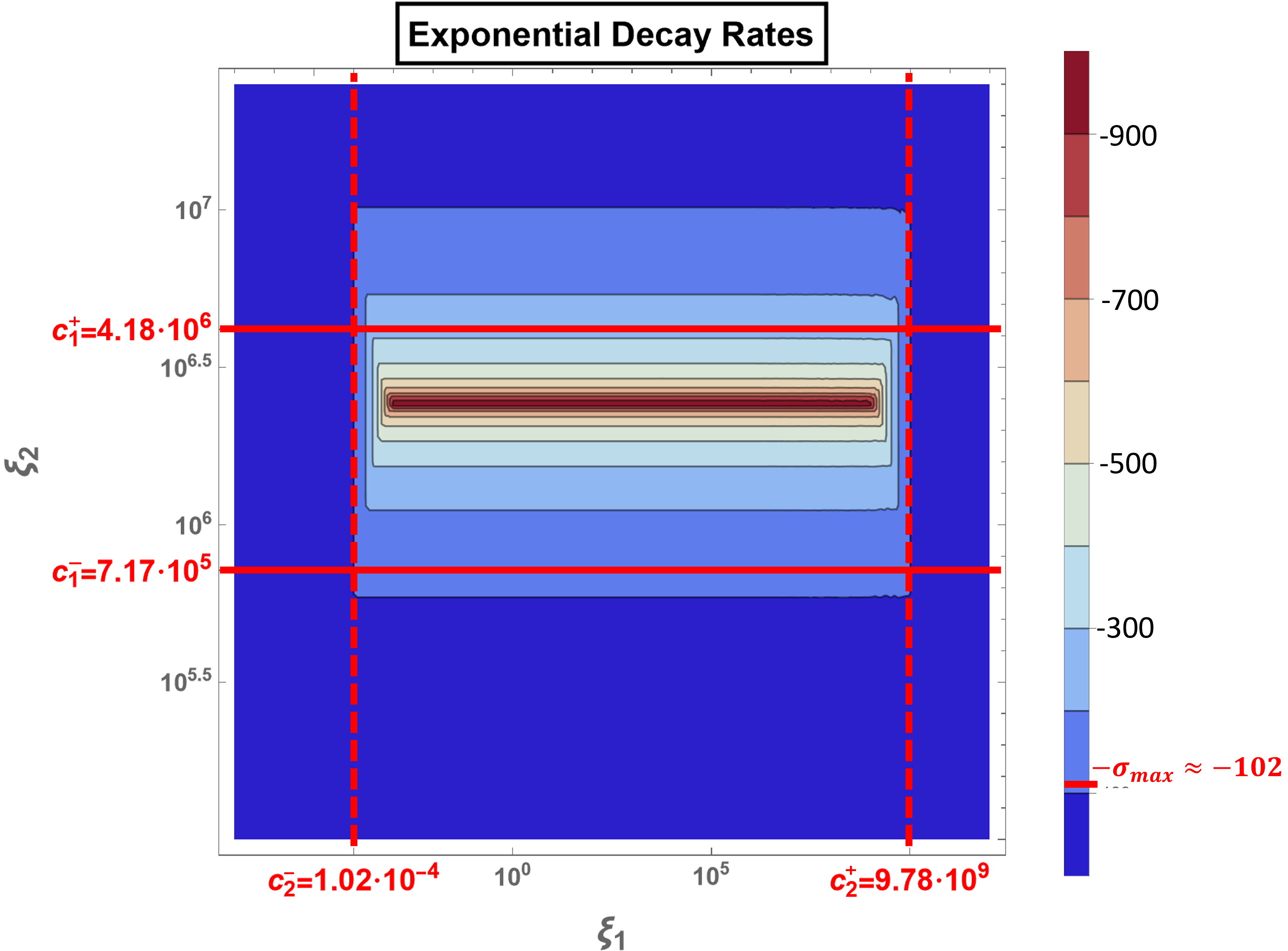}
\caption{
Decay rates for the ORFD approximation of the system \eqref{eq1} with \eqref{c1}, showing dependence on feedback amplifiers $\xi_1$ and $\xi_2$. Red lines mark the intervals $(c_1^-, c_1^+)$ and $(c_2^-,c_2^+)$, where the maximal decay rate is achieved. Lighter blue and brown areas indicate better decay rates, while darker blue regions show slower rates outside these intervals.}
	\label{DecayRates}
\end{figure}

}

\end{document}